\documentclass[preprint,review,12pt]{elsarticle}

\makeatletter
\def\ps@pprintTitle{%
  \let\@oddhead\@empty
  \let\@evenhead\@empty
  \def\@oddfoot{\reset@font\hfil\thepage\hfil}
  \let\@evenfoot\@oddfoot
}
\makeatother

\usepackage[utf8]{inputenc}
\usepackage{fullpage,times}
\usepackage{amsmath,amssymb}
\usepackage{graphicx}
\usepackage[ruled,vlined,linesnumbered,noend]{algorithm2e}
\usepackage{subfig}
\usepackage{caption}
\usepackage{url}
\usepackage[numbers]{natbib} 
\usepackage{placeins}
\usepackage{siunitx}
\usepackage{url}
\usepackage{xspace}
\usepackage[table,xcdraw]{xcolor}

\usepackage{hyperref}
\hypersetup{
colorlinks={true},
linkcolor={Blue3},
urlcolor={Blue3},
citecolor={Blue3},
pdfstartview={FitH},
pdfauthor={Author},
pdftitle={Title},
pdfsubject={Subject}
plainpages=false
}

\usepackage{comment}

\usepackage{enumitem}

\usepackage{booktabs} 
\usepackage{multirow}
\usepackage{multicol}
\usepackage{makecell}

\usepackage[toc,page]{appendix}

\usepackage{setspace}

\definecolor{wisconsin-red}{rgb}{0.6,0,0}
\definecolor{darkgreen}{rgb}{0.2,0.6,0.2}
\definecolor{maroon}{rgb}{0.5, 0.0, 0.0}
\definecolor{violet}{rgb}{0.75, 0.0, 1.0}
\definecolor{lightgray}{gray}{0.9}
\definecolor{navyblue}{rgb}{0.0, 0.0, 0.5}
\definecolor{darkmidnightblue}{rgb}{0.0, 0.2, 0.4}
\definecolor{midnightblue}{rgb}{0.0,0.4,0.85}
\definecolor{Gray}{gray}{0.75}
\definecolor{darkgreen}{rgb}{0,0.5,0}
\definecolor{apricot}{rgb}{0.98, 0.81, 0.69}

\newcolumntype{C}[1]{>{\centering\arraybackslash}p{#1}}
\newcolumntype{P}[1]{>{\raggedright\arraybackslash}p{#1}}
\newcolumntype{L}[1]{>{\raggedleft\arraybackslash}p{#1}}

\usepackage{upgreek}

\newcommand{\NeurADP}{\texttt{NeurADP}}

\newcommand{\MyopicILP}{\texttt{Myopic-ILP}}
\newcommand{\MyopicHeuristic}{\texttt{Myopic-Heuristic}}

\newcommand{\Time}{\mathcal{T}}

\newcommand{\DecisionSet}{\mathbf{A}}

\newcommand{\WorkersSet}{\mathcal{W}}
\newcommand{\OrdersSet}{\mathcal{O}}
\newcommand{\TimeIntervals}{\texttt{$\delta$}}

\newcommand{\CurrentTime}{\texttt{$t$}}
\newcommand{\StateNote}{\texttt{$S$}}

\newcommand{\ExogenousInformation}{W}

\newcommand{\Reward}{R}

\newcommand{\PostDecision}{\texttt{Post}}

\newcommand{\MaximumDelay}{\gamma}

\newcommand{\MValue}{\beta}
\newcommand{\IndividualWorker}{\texttt{$w$}}

\newcommand{\IndividualOrder}{\texttt{$o$}}
\newcommand{\WorkerHuman}{\text{human}}
\newcommand{\WorkerLocation}{\text{loc}}
\newcommand{\WorkerCapacity}{\text{cap}}
\newcommand{\WorkerBattery}{\text{bat}}
\newcommand{\WorkerOrders}{\text{ords}}
\newcommand{\OrderPickup}{\text{pickup}}
\newcommand{\OrderHuman}{\text{human}}
\newcommand{\OrderDeadline}{\text{dead}}
\newcommand{\OrderPost}{\texttt{Order-Post}}
\newcommand{\WorkerPost}{\texttt{Worker-Post}}

\newcommand{\SingleDecision}{a}

\begin{document}\sloppy
\setlength{\parindent}{2em}

\begin{frontmatter}
\title{Dynamic AGV Task Allocation in Intelligent Warehouses} 

\author[1]{Arash Dehghan}
\ead{arash.dehghan@torontomu.ca}

\author[1]{Mucahit Cevik\corref{cor1}%
\fnref{fn1}}
\ead{mcevik@torontomu.ca}

\author[2]{Merve Bodur}
\ead{merve.bodur@ed.ac.uk}

\cortext[cor1]{Corresponding author}
\fntext[fn1]{Toronto Metropolitan University, Toronto, ON, Canada}
\address[1]{Toronto Metropolitan University, Toronto, ON, Canada}
\address[2]{University of Edinburgh, Edinburgh, UK}

\begin{abstract}
This paper explores the integration of Automated Guided Vehicles (AGVs) in warehouse order picking, a crucial and cost-intensive aspect of warehouse operations. The booming AGV industry, accelerated by the COVID-19 pandemic, is witnessing widespread adoption due to its efficiency, reliability, and cost-effectiveness in automating warehouse tasks. This paper focuses on enhancing the picker-to-parts system, prevalent in small to medium-sized warehouses, through the strategic use of AGVs. We discuss the benefits and applications of AGVs in various warehouse tasks, highlighting their transformative potential in improving operational efficiency. We examine the deployment of AGVs by leading companies in the industry, showcasing their varied functionalities in warehouse management. Addressing the gap in research on optimizing operational performance in hybrid environments where humans and AGVs coexist, our study delves into a dynamic picker-to-parts warehouse scenario. We propose a novel approach Neural Approximate Dynamic Programming approach for coordinating a mixed team of human and AGV workers, aiming to maximize order throughput and operational efficiency. This involves innovative solutions for non-myopic decision making, order batching, and battery management. We also discuss the integration of advanced robotics technology in automating the complete order-picking process. Through a comprehensive numerical study, our work offers valuable insights for managing a heterogeneous workforce in a hybrid warehouse setting, contributing significantly to the field of warehouse automation and logistics.
\end{abstract}
\begin{keyword}
Warehouse management \sep Dynamic task allocation \sep Approximate dynamic programming \sep Deep learning
\end{keyword}
\end{frontmatter}
\vspace{-8pt}

\section{Introduction}
The Automated Guided Vehicle (AGV) industry is flourishing within the warehousing sector, with over 100 manufacturers competing in this growing market \citep{agvnetwork}. Financial investment in AGVs continues to soar, projected to hit US\$2.74 billion by 2023, propelled further by the COVID-19 pandemic which accelerated investments in this domain \citep{elias2020agv}. A significant 53\% of third-party logistics providers see warehouse automation as a prime opportunity for 2023, showcasing a strong trend towards embracing these technologies \citep{moveragv}. This pivot towards automation is not only a reaction to persistent labor shortages, but also a strategic step to stay competitive against major retail conglomerates \citep{beasley2023warehouse}. Retail behemoth Walmart anticipates that about 65\% of its stores will incorporate automation by 2026, emphasizing the critical role of AGV technology in shaping the retail sector's logistical future \citep{reuters2023amazon}. This collective transition indicates that warehouse automation, especially through AGV adoption, is swiftly becoming an industry norm, addressing modern logistics and supply chain challenges effectively.

The adoption of AGVs is driven by their significant advantages, which include improved efficiency that speeds up warehouse operations, as well as consistent and reliable workflow management which reduces human error and enhances predictability \citep{nextsmartship2022agv}. Financially, they provide notable direct labor cost reductions and indirect savings through improved accuracy and speed, while safety enhancements and reductions in workplace accidents bolster their value proposition further \citep{benevides2020agv, nextsmartship2022agv, strom2018agv}. Furthermore, the agility of AGVs in adapting to varied tasks, their scalability for business growth, their contribution to better space utilization, and the relative ease of integration into existing systems fortify their appeal \citep{benevides2020agv}. AGVs also enable a more strategic allocation of labor, allowing workers to focus on more complex or customer-centric tasks \citep{mosca2018agv, strom2018agv}. Subsequently, the breadth of AGV applications within warehouses continues to expand, encompassing tasks such as item delivery, stock replenishment, order picking, and loading and unloading, which are all essential to efficient warehouse operations \citep{mhi2023supplychain, parimi2017agv}. This wide-ranging functionality of AGVs showcases their transformative potential and their rising status as a staple in warehouse management and logistics.

Given the numerous advantages and the versatility in application, several companies have already made significant strides in integrating AGVs into their warehouse workflows. Amazon, a trailblazer in warehouse automation, has deployed over 100,000 robots across its fulfillment centers, while Seegrid has introduced robust pallet trucks and tow tractors which can handle massive loads of up to 8,000 and 10,000 pounds, respectively \citep{elias2020agv}. Additionally, Linde Material Handling caters to a broad range of customer needs, offering a diverse portfolio of AGVs designed for the multifaceted demands of daily warehouse operations \citep{lindeagv}, and Dematic has positioned its AGVs as essential components for material transport, bridging the gap between production areas and storage facilities, as well as executing various other logistical tasks \citep{dematicagv}. 

Our work investigates the role of AGVs in the order picking aspects of warehouse operations, which involve individual items being selected and collected from storage locations to fulfill customer orders. Addressing this aspect of warehouse operations is crucial as it is the most expensive and time-consuming. More specifically, the picking process accounts for up to 55\% of the total warehouse operating costs \citep{tompkins2010facilities}, and as much as 60\% of all labor activities within the warehouse \citep{de2007design}. Two common methods exist in manual order picking: picker-to-parts systems, where workers travel to item locations, collect them, and move to the next location; and parts-to-picker systems, which bring goods from storage to designated picking areas for worker selection. Our paper primarily focuses on the picker-to-parts systems as they are employed by a vast majority of small to medium-sized retail and logistics firms \citep{xu2020research}, and up to 90\% of warehouses in the grocery sector \citep{kuhn2013integrative}.

The integration of robotics into order picking operations presents a compelling opportunity for enhancements in operational efficiency. Manual picker-to-parts tasks are among the most laborious in warehouses, often resulting in musculoskeletal disorders, low back pain, and other physical ailments which can diminish the efficiency of the picking system \citep{vijayakumar2021literature}. Implementing robots not only addresses these human-centric issues but also offers a substantial reduction in costs related to human labor, such as insurance, salaries, and benefits \citep{benevides2020agv, strom2018agv}. Numerous studies have explored the integration of Autonomous Mobile Robots (AMRs) in the picking process to aid human pickers with the transportation of items within a hybrid environment \citep{azadeh2020dynamic, loffler2023human, pugliese2022amr, srinivas2022collaborative}. In such setups, humans are responsible for picking the items from shelves, while AMRs are utilized solely for transportation, thereby automating only a portion of the picking process. However, thanks to advancements in object picking technology, there are now companies which automate the entirety of the robot-to-parts process. For instance, Magazino's TORU \citep{LogistikKnowHow2023} is a sophisticated logistics robot that adeptly navigates to shelves to retrieve items directly or to pull cartons toward itself. When dealing with mixed cartons, its integrated 2D and 3D cameras employ 3D computer vision to scan the shelf contents. After matching the items with its database, it selectively extracts the targeted goods. These are then stored internally in the robot's adaptable compartments. Equipped with numerous sensors, TORU safely operates alongside humans, enabling the flexible automation of tasks that were previously done manually. This technology can reduce picking costs by up to 40\% when compared to traditional methods. Similarly, Fetch Robotics has introduced its Fetch and Freight robots \citep{wise2016fetch}. These machines navigate through ADA-compliant buildings and are equipped with arms capable of reaching down to retrieve items from the floor, ensuring efficient item recovery. Their comprehensive sensor array allows for object perception, navigation, and manipulation in dynamic settings.

As the concept of a hybrid order-picking environment --where robots and humans work together in a shared warehouse space-- is relatively new, few studies explored the optimization of operational performance in such settings. That is, the collaboration between humans and robots in warehouses remains under-explored, leaving several unanswered questions about how to most effectively pair orders with this joint workforce to enhance efficiency. Our research specifically focuses on this emerging paradigm.
We explore the coordination of a mixed team of humans and AGVs navigating the warehouse space, with the objective of intelligently pairing incoming orders to these operatives, taking into account the battery management of AGVs. Our specific modeling objective is to enhance operational efficiency and maximize order throughput. To that end, we develop a novel Markov Decision Process (MDP) model and devise a Neural Approximate Dynamic Programming (NeurADP) framework~\citep{shah2020neural} to handle efficient batching and assignment of incoming orders while concurrently managing the battery life of the robotic workers. While there has been some recent research in this specific area within the AGV literature, it has mostly concentrated on static, predictable scenarios, often employing basic heuristics or rule-based approaches for assigning incoming orders to human and robot teams, and has also neglected the charging aspects of the robots. A summary of contributions to the existing literature is provided as follows.
\begin{itemize}
    \item We formulate the order picking problem as an MDP to account for the uncertainty of stochastic order arrivals in a hybrid environment with human and AGV workers.
    Our model extends beyond previous research for this problem by being the first to incorporate charging stations within the hybrid warehouse setting and introduce battery management decision-making for AGVs.
    
    \item We implement a NeurADP solution methodology, advancing beyond traditional myopic- and heuristic-based solutions used in previous work. Our experimental results demonstrate that NeurADP significantly outperforms myopic and heuristic-based methods, both in terms of the quantity of orders processed, as well as in the efficiency achieved in executing picking tasks.
    
    \item We provide managerial insights related to the hybrid order picking setting as well as analysis on the impact of various factors such as the number and types of workers, worker speed, delay time allowance, worker capacity, and the availability of orders for both humans and AGVs, along with the incorporation of order deadlines.
\end{itemize}

The remainder of the paper is organized as follows. Section~\ref{literaturereview} provides a comprehensive review of the relevant AGV literature, better positioning our research within the existing body of work. Section~\ref{problemdescription} provides a formal description of the problem setting for our problem. In Section~\ref{solutionmethodology}, we describe the 
solution methodology. Details regarding the datasets and benchmark policies used in the experiments are provided in Section~\ref{experimentalsetup}. The results of the computational experiments are presented in Section~\ref{results}, followed by a conclusion in Section~\ref{conclusion} that summarizes the research findings and suggests avenues for future research.

\section{Literature Review} \label{literaturereview}
AGV control problems may be broken down into five core tasks: task allocation, in which the goal is to optimally assign a set of tasks to a set of AGVs, localization, where the goal is to locate the exact location on a map, path planning, which seeks to generate an obstacle-free path from two locations, motion planning, which requires real time modifications of a planned path according to dynamic obstacles, and vehicle management, which focuses on the management of vehicles battery, error, and maintenance statuses. Our research particularly relates to task allocation and vehicle management aspects. 

Battery management is an important part of AGV operations. While many studies highlight the importance of integrating battery management into AGV decision-making, emphasizing its significant influence on system performance \citep{LEANH20061,VIS2006677}, the literature still largely overlooks this area \citep{de2020automated}. \citet{kawakami2012battery} investigate the importance of battery management in AGV systems for reducing costs and improving efficiency. They specifically focus on Valve-regulated Lead-Acid batteries, commonly used in AGVs, highlighting the need for appropriate charging intervals to prevent battery deterioration and extend battery life. Similarly, \citet{kabir2017} examine the effects of different routing techniques for battery management on the performance of AGVs and analyze how routing to charging stations can impact the overall productivity of a manufacturing facility. \citet{KABIR2018225} study how adjusting the battery charging durations of AGVs can enhance manufacturing capacities in the short term, and \citet{de2020resource} introduce an advanced decentralized method for optimizing the integration of charging stations into the existing optimal tour routes of AGVs. Different from these studies, which focus primarily on the battery management aspects of AGV operations, our work extends the scope by integrating battery management considerations with order-picking task allocation in a hybrid warehouse setting.

Solution approaches for AGV task allocation may be broken down into two: optimization-based and market-based. In optimization-based methods, an algorithm searches for an optimal solution in a solution space which maximizes a profit or minimizes a cost using global information and considering all constraints. Whereas, in market-based methods, an economic principle is used to solve the task allocation problem. In this literature review, we focus on optimization-based methods given that they are most relevant to our work; a breakdown of further market-based solutions may be found in \citep{de2020automated}. Task allocation problems span many domains such as manufacturing, healthcare, and robotics and can be solved using a wide range of methodologies such as exact algorithms \citep{atay2006, coltin2010mobile, giordani2010distributed}, heuristics \citep{kmiecik2010task, li2017multi, liu2012centralized, sarkar2018scalable}, dynamic programming \citep{ralevic2012dynamic}, and many others. Several papers have looked at the incorporation of hybrid warehouse settings for order picking task allocation systems.

Recent studies have examined the use of AMRs in hybrid environments to assist human pickers, where humans pick items from shelves and AMRs handle transportation \citep{azadeh2020dynamic, loffler2023human, pugliese2022amr, srinivas2022collaborative}. On the other hand, innovations in object picking technology, such as Magazino's TORU robot and Fetch Robotics' Fetch and Freight robots, have enabled full automation of the picking process \citep{LogistikKnowHow2023, wise2016fetch}. The idea of a hybrid order-picking environment, combining human and robot collaboration in shared warehouse spaces, is relatively new, leading to a research gap in optimizing operational performance in these contexts. However, recent studies have begun to explore and address this emerging field. We highlight these relevant studies in Table~\ref{table:references}, which is comprised of seven indicators that provide information about the problem setting and solution methodology. These are: ``Solution Technique'', which describes the approach used to solve the problem, ``Non-Myopic'', which indicates whether a myopic solution technique is employed, ``No Prior Knowledge'', which indicates whether orders are not known at the beginning of each work day, ``Deadline'', which indicates whether specific deadlines are set for fulfilling incoming orders, ``Shared Area'', which indicates whether human and AGVs workers share a workspace in the hybrid setting, ``Batching'', which indicates whether orders are able to be batched together, and finally ``Charging'', which indicates whether charging decisions and battery management is considered for AGVs.

\setlength{\tabcolsep}{2.5pt} 
\renewcommand{\arraystretch}{1.3} 
\begin{table}[!ht]
\centering
\caption{Summary of relevant studies.}\label{table:references}
\resizebox{0.99\textwidth}{!}{
\begin{tabular}{P{0.29\textwidth}C{0.13\textwidth}C{0.13\textwidth}C{0.13\textwidth}C{0.13\textwidth}C{0.13\textwidth}C{0.13\textwidth}C{0.13\textwidth}} 
\toprule
\textbf{Study} & \textbf{Solution Technique} & \textbf{Non-Myopic} & \textbf{No Prior Knowledge} & \textbf{Deadline} & \textbf{Shared Area} & \textbf{Batching} & \textbf{Charging} \\ 
\midrule
\citet{sgarbossa2020robot} & Heuristic  & & \checkmark & & & & \\
\citet{zhang2021evaluation} & Heuristic & &  \checkmark & & \checkmark & \checkmark & \\
\citet{kauke2022mobile} & Heuristic & & \checkmark & &    \checkmark & \checkmark & \\
\citet{winkelhaus2022hybrid} & Heuristic & & & & \checkmark & \checkmark & \\
\citet{zhang2023ergonomic} & Heuristic & & & & \checkmark & \checkmark & \\ 
\midrule
\textbf{Our Work} & NeurADP & \checkmark & \checkmark & \checkmark & \checkmark & \checkmark & \checkmark \\ 
\bottomrule
\end{tabular}
}
\end{table}

\citet{sgarbossa2020robot} introduce a robotic picker designed for pallet retrieval and formulate a strategy for allocating products between two distinct warehouse areas, which are designated respectively for human workers and robots. This approach involves a dual-objective optimization model aimed at reducing the labor intensity for human workers while simultaneously enhancing the uniformity of product categories assigned to each zone. To achieve these objectives, the authors employ the non-dominated sorting genetic algorithm to balance the minimization of human workload and the maximization of product category similarity within the designated zones. \citet{zhang2021evaluation} develop a simulation model to assess the energy expenditure of human pickers in a collaborative environment with picking robots, analyzing the operational costs, efficiency, and ergonomic impact. The model defines distinct roles for human pickers and robots based on a set of assignment rules that are contingent on different item classes. Under these rules, items belonging to specific classes are allocated either to human pickers or robots, depending on the nature of the item and the predefined criteria. The study's scenarios are then categorized and evaluated based on varying combinations of these assignment rules, leading to different distributions of tasks between humans and robots. \citet{kauke2022mobile} investigate the effects of aisle widths and layout variations on human-robot interaction in order picking systems, focusing on enhancing performance efficiency. Their study contrasts zoning strategies with traditional order picking systems, emphasizing the complexity and coordination demands of hybrid systems involving both humans and robots. Their approach uses a heuristic to determine picking routes, and assigning tasks to either humans or robots without considering individual traits. \citet{winkelhaus2022hybrid} present a simulation model for evaluating the performance of hybrid order picking systems, incorporating variables such as picker blocking. In their model, order assignments are based on predetermined workloads for each team, categorized by item classes (e.g., A, B, or C) and their turnover rates. Each item is pre-assigned to a specific team according to its class. When an item is required for an order, the designated team member, either a human or a robot, is responsible for picking it. \citet{zhang2023ergonomic} propose an agent-based simulation model to explore how hybrid order picking systems can reduce the daily workload of human pickers. Their model operates under the assumption that all customer orders are known before each shift, eliminating idle time due to late-arriving orders. Orders are assigned following a ``first-come-first-served'' principle. The assignment rules stipulate that human pickers handle `A' items, while robots are responsible for `B' and `C' items.

Our work explores the combined use of humans and AGVs in a warehouse to process orders. It extends previous research by not only pairing orders with workers but also by managing AGV charging. Unlike previous heuristic-based methods, we employ a non-myopic NeurADP approach, which is shown to be effective in such complex scenarios~\citep{shah2020neural}. 
Additionally, our model incorporates order deadlines and our empirical study includes a detailed sensitivity analysis on various key parameters such as worker types, speeds, delay allowances, capacity, and order availability. As such, our paper provides new managerial insights for effectively operating in this hybrid warehouse setting.

\section{Problem Description and Formulation} \label{problemdescription}
In our study, we introduce a dynamic picker-to-parts model tailored for a hybrid warehouse environment, integrating both human and AGV workers. This model aims to efficiently allocate workers to incoming order batches and determine optimal charging strategies for AGVs, including the timing and duration of charging sessions. It functions over a 24-hour decision horizon, adapting to the variable demand patterns of orders in a grid-layout warehouse equipped with strategically placed charging stations and a designated drop-off area for picked items. Orders in our system are generated stochastically, each with its own delivery deadline influenced by its arrival time. Upon assignment to a worker, an order becomes an incorporated part of the system, with a guarantee to meet its designated drop-off deadline. Moreover, the model takes into account a predetermined group of heterogeneous workers available during the planning horizon, considering their capacity constraints and, in the case of AGVs, their battery levels.

Our model incorporates several key problem specifications related to the hybrid order-picking problem setting. First, workers may be matched with multiple orders per time-step, wherein each order is associated with its own pick-up location. A worker is then expected to traverse the warehouse to pick up the orders at each of their respective locations and deliver them to the designated drop-off area. Workers maintain a queue of their assigned orders to track which ones need to be picked up and which have been collected, before they are delivered to the drop-off area. The queue is dynamically rearranged whenever a new order or batch is assigned to a worker, optimizing their route within the warehouse for order collection and return to the drop-off area. Once an order is placed in a worker's queue, it cannot be transferred to another worker's queue. Unmatched orders that surpass their arrival period are removed from the system, reflecting the expectation of customers for timely confirmation of their requests. However, new orders arriving in subsequent time steps can be assigned to workers and added to their queues, provided several considerations are taken into account.

First, each worker has a designated capacity corresponding to the size of their storage bin, which is used to hold orders. Subsequently, a batch of orders is only eligible to be matched with a worker if the additional capacity assigned to them does not exceed their bin capacity. Furthermore, a batch of orders may only be assigned to a worker and added to their queue if by adding the batch of orders, the worker is still able to pick up all remaining orders and drop off all orders at the drop-off area prior to each orders deadline. For AGVs, battery levels are additionally a crucial factor in matching them with order batches. Specifically, an AGV is only assigned a batch if it can complete all orders in that batch, as well as those already in its queue, and still reach the nearest charging station before its battery depletes completely. This policy ensures AGVs maintain sufficient charge throughout the operational period. Additionally, it is assumed that AGVs cannot charge while serving orders. Once a decision to charge is made, the AGV will charge uninterrupted until the start of the next decision epoch. 

The main goal of our model is to maximize the fulfillment of online orders within the decision horizon. We factor in uncertainties of future order arrivals and the downstream effects of current choices, including those related to charging decisions. To manage these complex decisions, we develop an MDP model and incorporated a NeurADP solution framework. This approach facilitates efficient real-time decision-making, even amid uncertainty. We present the MDP model components in Table~\ref{table:mdpmodel}.

\setlength{\tabcolsep}{4.5pt}
\renewcommand{\arraystretch}{1.25}
\begin{table}[ht]
\centering
\caption{MDP model components}
\label{table:mdpmodel}
\resizebox{0.905\textwidth}{!}{
\begin{tabular}{|l|l|}
\hline
\textbf{Component} & \textbf{Notation/Description} \\ \hline
Decision epochs & $ \CurrentTime \in \Time = \{0, ..., T\} $ \\ \hline
System state = (Workers, Orders) & $ \StateNote_\CurrentTime = (\WorkersSet_\CurrentTime,\OrdersSet_\CurrentTime) $ \\ \hline
\quad Worker attributes & \quad $ \IndividualWorker = (\IndividualWorker_{\WorkerLocation}, \IndividualWorker_{\WorkerHuman}, \IndividualWorker_{\WorkerCapacity}, \IndividualWorker_{\WorkerBattery}, \IndividualWorker_{\WorkerOrders}) \in \WorkersSet $ \\ \hline
\quad Order attributes & \quad $ \IndividualOrder = (\IndividualOrder_{\OrderPickup}, \IndividualOrder_{\OrderHuman}, \IndividualOrder_{\OrderDeadline}) \in \mathcal{O} $\\ \hline
Decisions  & \( \SingleDecision_\CurrentTime \) = (Assign order batch, Charge AGV, Null) \\ \hline
Feasible decisions & \( \SingleDecision_\CurrentTime \in 
\DecisionSet_\CurrentTime(\StateNote_\CurrentTime) \) \\ \hline
Immediate reward & \( \Reward_\CurrentTime(\SingleDecision_\CurrentTime) \) = $\MValue~\cdot $ Orders fulfilled - Time until drop-off \hfill \\ \hline
Post-decision state & \( \StateNote_{\CurrentTime}^{\PostDecision} = \texttt{statepost}(\StateNote_\CurrentTime,\SingleDecision_\CurrentTime) = (\StateNote_{\CurrentTime}^{\WorkerPost}, \StateNote_{\CurrentTime}^{\OrderPost} =  \emptyset) \) \\ \hline
Exogenous uncertainty information & $\ExogenousInformation_{\CurrentTime+1}$:  orders arriving between \( \CurrentTime \) and \( \CurrentTime+1 \) \\ \hline
State transition  & \( \StateNote_{\CurrentTime+1} = \texttt{statenext}(\StateNote_\CurrentTime^{\PostDecision}, \ExogenousInformation_{\CurrentTime+1}) = (\StateNote_{\CurrentTime}^{\WorkerPost},\ExogenousInformation_{\CurrentTime+1}) \) \\ \hline
Value function & \( V_\CurrentTime(\StateNote_\CurrentTime) \) \\ \hline
Post-decision value function & \( V_\CurrentTime^\PostDecision(\StateNote_\CurrentTime^\WorkerPost) = \mathbb{E}_{\ExogenousInformation_{\CurrentTime+1}} [V_{\CurrentTime+1}(\StateNote_{\CurrentTime+1}) | \StateNote_\CurrentTime^\WorkerPost] \) \\ \hline
Bellman optimality equation & $V_\CurrentTime(\StateNote_\CurrentTime) = \max \{ \Reward_\CurrentTime(\SingleDecision_\CurrentTime) + V_\CurrentTime^\PostDecision(\StateNote_\CurrentTime^\WorkerPost) : \SingleDecision_\CurrentTime \in \DecisionSet_\CurrentTime(\StateNote_\CurrentTime) \}$  \\ \hline
MDP formulation & $\max V_0(\StateNote_0)$ \\ \hline
\end{tabular}
}
\end{table}

To begin, the planning horizon $\Time$ is divided into discrete time intervals of $\TimeIntervals$ 
(e.g., five minutes), such that decisions are made at the beginning of each interval, while exogenous information, denoted by $\ExogenousInformation$, is observed continuously throughout. The state of the system $\StateNote_\CurrentTime$ is characterized by the attributes of workers and orders, defined by $\WorkersSet$ and $\OrdersSet$, respectively. An individual worker's state is captured by a five-dimensional attribute vector $\IndividualWorker$. This vector includes $\IndividualWorker_\WorkerLocation$, indicating their position in the warehouse, and $\IndividualWorker_\WorkerHuman$, a binary attribute distinguishing between human workers and AGVs. $\IndividualWorker_\WorkerCapacity$ reflects the worker's current order capacity, $\IndividualWorker_\WorkerBattery$ shows their battery level (applicable only to AGVs), and $\IndividualWorker_\WorkerOrders$ details the worker’s task queue, including assigned orders. This queue, optimized for minimal travel time, is updated with each new batch of orders assigned to the worker. Moreover, the state of each order is represented by a three-dimensional vector $\IndividualOrder$. This includes $\IndividualOrder_\OrderPickup$, specifying the order's storage location in the warehouse, $\IndividualOrder_\OrderHuman$, a binary value indicating if the order requires handling exclusively by humans, and $\IndividualOrder_\OrderDeadline$, the drop-off deadline, calculated as $\CurrentTime + \MaximumDelay_\IndividualOrder$. Here, $\gamma_\IndividualOrder$ denotes the permissible delay time for the order $\IndividualOrder$.

The set of feasible decisions may be defined by $\SingleDecision_\CurrentTime \in \DecisionSet_\CurrentTime(\StateNote_\CurrentTime)$. An individual action $\SingleDecision_\CurrentTime$ for a worker may encompass one of the following: allocating a batch of orders, designating a charging task, or assigning a null action. The assignment of order batches depends on the worker's capacity and ability to deliver both new and existing orders before their deadlines. For AGVs, this also includes the need to have sufficient battery to complete deliveries and potentially reach a charging station complete depletion. The charging action, specific to AGVs, is only feasible if the AGV is not actively serving orders. When assigned a charging action, the worker heads to the nearest charging station in the warehouse. If they arrive at the charging station before the next decision epoch, they use the remaining time until that epoch to recharge their battery. Conversely, if the worker doesn't reach the charging station by the next decision epoch, a new decision will be made for them in the subsequent epoch, based on their updated state. Finally, the null action implies that the worker continues their current activity, whether that is progressing with order pick-up and delivery or remaining idle in the absence of new orders. Furthermore, the immediate reward for an action $\SingleDecision_\CurrentTime$ is determined by multiplying the total orders served by $\MValue$ and then subtracting the time from $\CurrentTime$ until the worker can deliver their assigned orders to the drop-off area. $\MValue$ is introduced to prioritize maximizing the number of orders completed per interval, ensuring this aspect outweighs the time factor. When actions serve an equal number of orders, preference is given to the one that minimizes task completion time, thereby freeing workers sooner for new order batches. This aspect is calculated by multiplying the worker's maximum capacity length by $\MaximumDelay$.

\begin{figure}[htb]
\centering
\includegraphics[width=\textwidth]{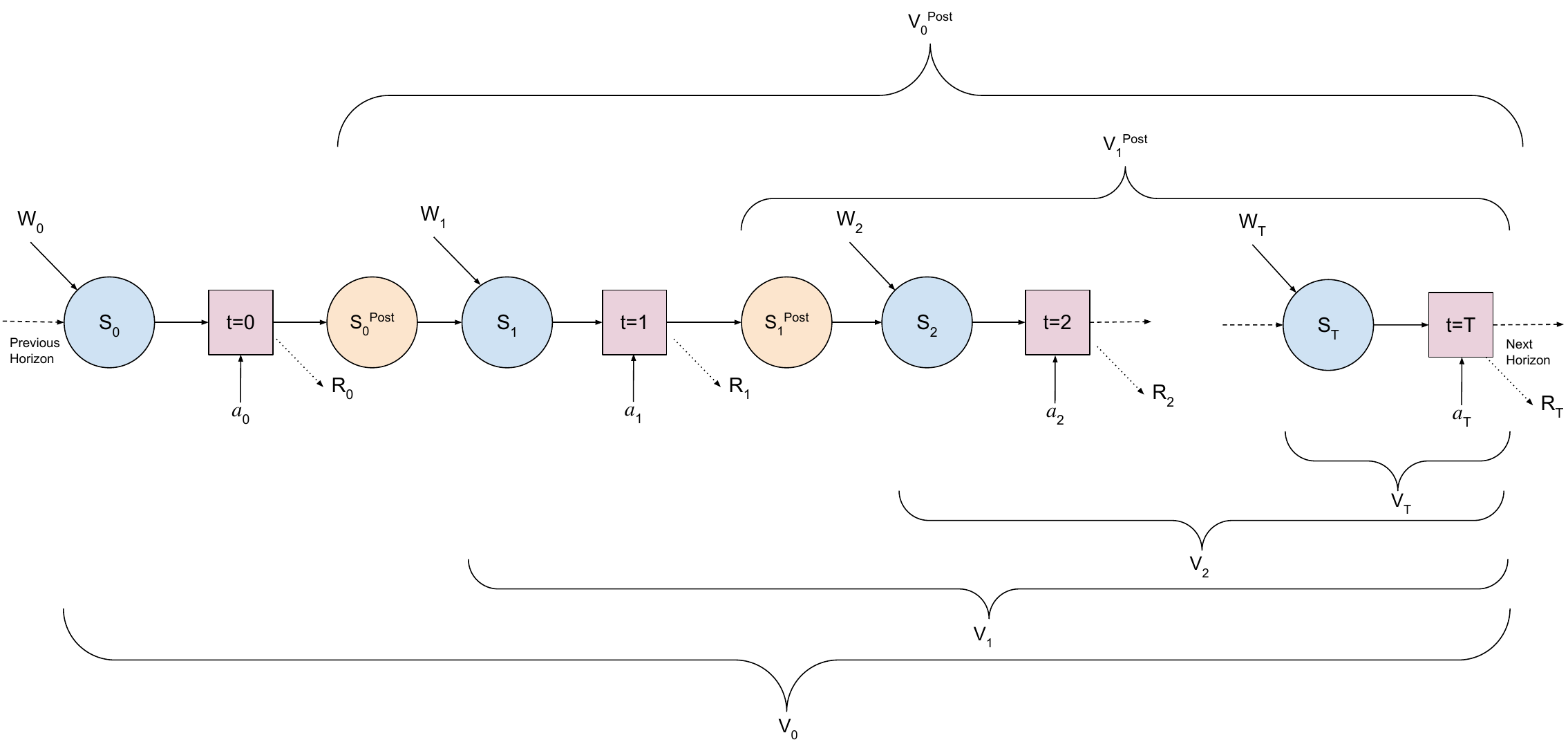}
\caption{System evolution for one planning horizon}
\label{fig:adp_evolution}
\end{figure}

In MDP models, the system evolution involves the transition from an initial state $\StateNote_0$, via a set of actions $\SingleDecision_0$, to a subsequent state $\StateNote_1$. This recursive evolution is continued until the final decision epoch of the horizon at time $T$. However, as illustrated in Figure~\ref{fig:adp_evolution}, the system evolution in ADP is more granularly defined to explicitly denote pre-decision and post-decision states, as well as the arrival of exogenous information. Given an initial state $\StateNote_0$, and the arrival of exogenous information $\ExogenousInformation_0$, a set of actions $\SingleDecision_0$ may be taken at time $\CurrentTime=0$, corresponding to subsequent rewards of $\Reward_0$. The system then progresses to a post-decision state, denoted as $\StateNote_0^\PostDecision=\texttt{statepost}(\StateNote_0,\SingleDecision_0)$. This represents the system's state after implementing actions on $\StateNote_0$, but before the arrival of new exogenous information in the next time step. The subsequent state, $\StateNote_1$, emerges through the receipt of exogenous uncertainty $\ExogenousInformation_1$ and the state transition function $\texttt{statenext}(\StateNote_0^\PostDecision,\ExogenousInformation_1)$. This recursive process again repeats up to the final decision epoch at time $T$. In our order picking problem, given that orders exit the system at the end of each decision epoch, the post-decision state may be represented as $\StateNote^\PostDecision_\CurrentTime=(\StateNote_{\CurrentTime}^{\WorkerPost}, \StateNote_{\CurrentTime}^{\OrderPost} = \emptyset)$. The state transition is then expressed by $\StateNote_{\CurrentTime + 1}=(\StateNote^\WorkerPost_\CurrentTime,\ExogenousInformation_{\CurrentTime + 1})$, where $\ExogenousInformation_{\CurrentTime + 1}$ denotes the arrival of orders between $\CurrentTime$ and $\CurrentTime + 1$.

Given that $V_\CurrentTime(\StateNote_\CurrentTime)$ denotes the value of being in a state $\StateNote_\CurrentTime$ at decision epoch $\CurrentTime$, we may define the Bellman optimality equation as such:
\begin{equation}
    V_\CurrentTime(\StateNote_\CurrentTime) = \max \{ \Reward_\CurrentTime(\SingleDecision_\CurrentTime) + \mathbb{E}_{\ExogenousInformation_{\CurrentTime+1}} [V_{\CurrentTime+1}(\StateNote_{\CurrentTime+1}) | \StateNote_\CurrentTime, \SingleDecision_\CurrentTime, \ExogenousInformation_{\CurrentTime + 1}] : \SingleDecision_\CurrentTime \in \DecisionSet_\CurrentTime(\StateNote_\CurrentTime) \}.
    \label{eq:Bellman}
\end{equation}
Incorporating the post-decision state, we may break down and rewrite the Bellman equation as follows:
\begin{subequations}
\label{m:BellmanBreakup}
    \begin{align}
    V_\CurrentTime(\StateNote_\CurrentTime) = \max \{ \Reward_\CurrentTime(\SingleDecision_\CurrentTime) + V_\CurrentTime^\PostDecision(\StateNote_\CurrentTime^\WorkerPost) : \SingleDecision_\CurrentTime \in \DecisionSet_\CurrentTime(\StateNote_\CurrentTime) \} \label{eq:BellOne} \\
    V_\CurrentTime^\PostDecision(\StateNote_\CurrentTime^\WorkerPost) = \mathbb{E}_{\ExogenousInformation_{\CurrentTime+1}} [V_{\CurrentTime+1}(\StateNote_{\CurrentTime+1}) | \StateNote_\CurrentTime^\WorkerPost, \ExogenousInformation_{\CurrentTime + 1}]  \label{eq:BellTwo}
    \end{align}
\end{subequations}
To simplify calculations and avoid the computational burden of examining each possible outcome for future values, Equation~\eqref{eq:BellOne} formulates a deterministic optimality equation based on the post-decision state, while Equation~\eqref{eq:BellTwo} defines the value of the post-decision state as the expected total of future rewards. Given the impracticality of calculating this expected value, an approximate value of the post-decision state value function is taken instead. This approach facilitates a more manageable and efficient estimation of future rewards, streamlining the optimization process in large-scale problems such as the one in this paper. To determine the optimal policy which maximizes the total expected reward, we must maximize the initial value function $V_0$, given the initial state $\StateNote_0$. As such, the objective of our MDP model can be defined as $\max V_0(\StateNote_0)$.
\section{Solution Methodology} \label{solutionmethodology}

\newcommand{\powerset}{\mathcal{P}}
We adapt the NeurADP framework \citep{shah2020neural} to our problem. NeurADP is a novel ADP-based algorithm designed for large-scale decision-making problems. It employs neural network value function approximations and utilizes deep reinforcement learning techniques for improved stability and efficiency. Furthermore, its ability to learn from integer programming-based assignments to manage complex combinatorial challenges makes it useful for large-scale problems. 

In what follows, we explain the details of the adaptation of the general NeurADP framework to the dynamic AGV task allocation problem, to derive high-quality approximate solutions to our proposed MDP model. We start with the building blocks and then provide the overall NeurADP algorithm. Given system state $S_t = (W_t,O_t)$, we mainly follow the below steps.   
\begin{enumerate}
\item {\it Enumerate the set of feasible order batches for workers:} Let $\Gamma_t: W_t \rightarrow \powerset^2(O_t)$ where $ \powerset$ denotes the power set and $ \powerset^2(\cdot) =  \powerset ( \powerset (\cdot))$. Given a worker $w \in W_t$, this function returns the collection of order batches (created from the orders in $O_t$) feasible to be assigned to the worker. That is, $\Gamma_t(w)$ is a set with each element $\mathcal{G} \in \Gamma_t(w)$ corresponding to a batch of orders (i.e., a subset of $O_t$) such that the worker of attributes $w$ can handle all of its currently assigned tasks combined with all the orders in $\mathcal{G}$ in a feasible manner. 

Algorithm~\ref{algo:Matchingfeasibility} defines the enumeration of the matching feasibility set $\Gamma_\CurrentTime(\IndividualWorker)$. The algorithm accepts as input the attribute state of a worker, $\IndividualWorker$, and the set of orders, $\OrdersSet_\CurrentTime$. It begins by initializing an empty set for the feasible matchings. Next, provided that the worker has available capacity, the set of potential feasible batchings $\powerset(\OrdersSet_\CurrentTime)$ is iterated. If the worker is an AGV and there exists an order within potential batching which may only be handled by humans, said batching is ignored. Otherwise, the feasibility of matching potential batching $\mathcal{G}$ to worker $\IndividualWorker$ is derived through Algorithm~\ref{algo:feasibility}. The algorithm begins by combining the set of potential batch orders, $\mathcal{G}$, with the orders already assigned to a worker, $\IndividualWorker_\WorkerOrders$. This results in a new set of assigned orders, $\IndividualWorker'_\WorkerOrders$. Additionally, $\IndividualWorker'_\texttt{pickups}$ is initialized to denote the new remaining set of pick-up locations the worker must visit prior to dropping orders off at the drop-off area. The algorithm then proceeds to iterate through all possible permutations of routes for collecting each order in the set $\IndividualWorker'_\texttt{pickups}$ and dropping them off. If a path exists where all orders are delivered before their respective deadlines, and the worker is either human or an AGV with adequate battery life to complete the route and recharge afterwards, then the algorithm confirms the feasibility of the path. Otherwise, it returns false. If Algorithm~\ref{algo:feasibility} returns a value of true, then the potential batching $\mathcal{G}$ in Algorithm~\ref{algo:Matchingfeasibility} is added to the matching feasibility set $\Gamma_\CurrentTime(\IndividualWorker)$, otherwise, the algorithm iterates to the next potential batching. Algorithm~\ref{algo:Matchingfeasibility} ends by returning the matching feasibility set $\Gamma_\CurrentTime(\IndividualWorker)$.


\begin{algorithm}[!ht]
\setstretch{1.23}
  \small
  \SetAlgoLined
    \SetKwInOut{Input}{Input}
    \SetKwInOut{Output}{Output}
    \Input{$w$, $O_t$}
    \Output{Set of feasible matchings}
    \textbf{MatchingFeasibility}($w,O_t$):\\
    \-\hspace{0.4cm} Initialize $\Gamma_t(w)$ as an empty set\;
    \-\hspace{0.4cm} \uIf{$w_{\texttt{cap}} < \texttt{maxcap}_w$ }{
    \-\hspace{0.4cm}     \For{$\mathcal{G} \in \powerset(O_t)$}{
    \-\hspace{0.4cm}     \uIf{$w_{\texttt{human}}== \texttt{false}$ and $o_{\texttt{human}}== \texttt{true}$ for any $o\in \mathcal{G}$}{
    \-\hspace{0.4cm}     continue
    }
    \-\hspace{0.4cm}     \uIf{isFeasible($w, \mathcal{G}$) $==true$}{
    \-\hspace{0.4cm}     $\Gamma_t(w) \leftarrow \Gamma_t(w) \cup \{ \mathcal{G}\}$
        }
        }
    }
    \-\hspace{0.4cm} \Return{$\Gamma_t(w)$}
     \caption{Matching feasibility}
     \label{algo:Matchingfeasibility}
\end{algorithm}

\begin{algorithm}[!ht]
\setstretch{1.23}
  \small
  \SetAlgoLined
\SetKwInOut{Input}{Input}
\SetKwInOut{Output}{Output}
\Input{$w$, $\mathcal{G}$}
\Output{Boolean value indicating if the path is feasible}
\textbf{isFeasible}($w, \mathcal{G}$):\\
\-\hspace{0.4cm} initialize $w'_{\texttt{ords}} \leftarrow w_{\text{ords}} \cup \mathcal{G}$\;
\-\hspace{0.4cm} initialize $w'_{\texttt{pickups}}$ as set of pick-up locations for orders in $w'_{\texttt{ords}}$ needing to be picked up\;
\-\hspace{0.4cm} \For{each path in permutation of $w'_{\texttt{pickups}}$}{
\-\hspace{0.4cm}     \uIf{path$_{\texttt{dropoff-time}} \leq o_{\texttt{dead}}$ for all $o \in w'_{\texttt{ords}}$}{
\-\hspace{0.4cm}         \uIf{$w_{\texttt{human}}$ == true or $w_{\texttt{bat}}$ is sufficient for battery required for action}{
\-\hspace{0.4cm}             \Return{true}
        }
    }
}
\-\hspace{0.4cm} \Return{false}
\caption{Path feasibility}
\label{algo:feasibility}
\end{algorithm}

\item {\it Define decision variables:} For worker $w \in W_t$, in addition to the possible order batch assignment, we have possible actions of going to the nearest charging station (for AGVs) and taking the null action, i.e., just continue with the previously assigned tasks. In that regard, we define the following binary decision variables: $x_{w \leftrightarrow \mathcal{G}}$ which takes the value of 1 if the order batch $\mathcal{G}$ is assigned, and 0 otherwise; $y_w$ which takes the value of 1 if the null action at the current state is instructed, and 0 otherwise; and $z_w$ which takes the value of 1 if the charging action is taken, and 0 otherwise.
\item {\it Create task allocation model:} Using the above-defined decision variables, we can obtain a binary programming representation of the feasible set of the Bellman optimality equation, previously stated as  
\begin{align}
V_t(S_t) =  \max \ & R_t(a_t) + V_t^\texttt{Post}(S_t^\texttt{Worker-Post}) \\
\text{s.t.} \ & a_t \in A_t(S_t),
\end{align}
which yields the following task allocation model: 
\begin{alignat}{2}
\max \ & \ R_t(x,y,z) + V_t^\texttt{Post}(W_t,x,y,z) \nonumber \\
\text{s.t.} \ & \sum_{\mathcal{G} \in \Gamma_t(w)} x_{w \leftrightarrow \mathcal{G}} + y_w = 1 && \forall w \in W_t : w_\text{human} = \text{True} \label{m:IPhumanactions}
 \\
& \sum_{\mathcal{G} \in \Gamma_t(w)} x_{w \leftrightarrow \mathcal{G}} + y_w + z_w = 1 \qquad && \forall w \in W_t : w_\text{human} = \text{False}  \label{m:IPhumanactions2} \\
& \sum_{w \in W_t} \sum_{\mathcal{G} \in \Gamma_t(w) : o \in \mathcal{G}} x_{w \leftrightarrow \mathcal{G}} \leq 1 && \forall o \in O_t  \label{m:orders}\\
& z_w = 0 && \forall w \in W_t : w_\text{human} = \text{True} \label{m:IPhumanactions3}\\
& \ x,y,z \ \text{binary} \label{m:variableDefinitions}
\end{alignat}
The model assigns one feasible action to each worker (namely order batch, null action, or charging; the last being eligible only for AGVs), ensuring that each order is assigned to at most one worker. With a slight abuse of notation, we parametrized the immediate and expected future reward of these actions by the binary decision variables as well as the workers' state vector for the latter. 

\item {\it Approximate the objective function:} We use a linear approximation for the objective function of the task allocation model:
\begin{align}
\max \ \sum_{w \in W_t} \sum_{\mathcal{G} \in \Gamma_t(w)} \alpha^\texttt{x}_{w \leftrightarrow \mathcal{G}} x_{w \leftrightarrow \mathcal{G}} \ + \sum_{w \in W_t} (\alpha^\texttt{y}_w y_w + \alpha^\texttt{z}_w z_w)
\label{m:IPobj}
\end{align}
where the coefficients 
$(\alpha^\texttt{x},\alpha^\texttt{y},\alpha^\texttt{z})$ are predicted via a neural network (NN). Note that since the charging action is not relevant to human workers, we can ignore those decisions in the objective, i.e., treat $\alpha^\texttt{z}_w = 0$ as fixed for any human worker index $w$; just for the ease of presentation we use the full coefficient vector.  As a result, we solve the \emph{approximated task allocation model} given by \eqref{m:IPhumanactions}-\eqref{m:IPobj} for decision making. 

We note that having this separable objective form, the NN is only providing estimates for the post-decision value function \emph{per worker}, not for the joint value of all the workers. More specifically, instead of approximating $V_t^\texttt{Post}(S_t^\texttt{Worker-Post})$, it helps with the approximation of the form 
$$\sum_{w \in \mathcal{W}} V^\approx_t \left(\{x_{w \leftrightarrow \mathcal{G}}\}_{\mathcal{G} \in \Gamma_t(w)},y_w,z_w,W_t \right).$$
The function $V^\approx_t$ basically takes the input of the action selected for a fixed worker $w$ (since only one the corresponding $x$, $,y$, and $z$ variables should be selected in the task allocation model) along with the current (i.e., pre-decision) state of all the workers. As a result, it has access to the post-decision state of only worker $w$ -- that is, $S_{tw}^\texttt{Worker-Post}$-- and uses the pre-decision state of the other workers -- that is, $W_t$ from $S_t$-- as auxiliary information to potentially improve the prediction for worker $w$. 

\item {\it Learn objective coefficient prediction:} We begin by sampling an experience, containing the state of workers $\WorkersSet_\CurrentTime$, the associated feasible action set, and the post-decision state of the workers from the previous time step $\StateNote^\PostDecision_{\CurrentTime - 1}$. The experience is then evaluated by applying each feasible action to the state of workers and scoring the post-decision state reached. This scoring for each post-decision state is carried out using a target neural network. The task allocation model is then employed to determine the actions of each worker. Subsequently, the post-decision states of the workers $\StateNote^\PostDecision_{\CurrentTime - 1}$ from the previous time step are updated through gradient descent. This update utilizes the supervised target scores obtained from the task allocation model.
\end{enumerate}

\begin{algorithm}[!ht]
\setstretch{1.23}
\small
\SetAlgoLined
\SetKwInOut{Input}{Input}
\SetKwInOut{Output}{Output}
\SetKwFor{For}{for}{do}{end for}
\SetKwRepeat{Do}{do}{while}
\SetKw{KwTo}{in}
\SetKwIF{If}{ElseIf}{Else}{if}{then}{else if}{else}{endif}
\SetKwFor{While}{while}{do}{endw}
\Input{Initial state of the system, $S_0$, and Neural Network (NN)}
\Output{Updated NN after simulation}
\BlankLine
\Do{a stopping criterion is not met}{
    \tcp*[h]{Simulate the system for one planning horizon}\;
    Initialize workers' and orders' states, $S_0$\;
    \For{$t=0,\hdots, T$}{
        Sample the new orders, $W_t$\;
        Enumerate the set of feasible order batches for workers, $\Gamma_t(\cdot)$; store them as an experience\;
        Obtain the decision-making model objective coefficients from the currently trained version of the NN, $(\alpha^{\texttt{x}}, \alpha^{\texttt{y}}, \alpha^{\texttt{z}})$\;
        Create the approximated task allocation model (using the predicted coefficients) and obtain its optimal solution\;
        (Optional) Sample from previously collected experiences and update the NN\;
        Update the system state implementing the obtained task allocation solution, $S_{t+1}$\;
    }
}
\caption{NeurADP training for AGV task allocation}
\label{algo:neural_network_training}
\end{algorithm}

Algorithm~\ref{algo:neural_network_training} defines the NeurADP training algorithm for the AGV task allocation problem. The algorithm takes as input the initial state of the system $\StateNote_0$, as well as an initialized neural network function. It begins by initializing the states of the workers, as well as the initial set of incoming orders, for the initial state $\StateNote_0$. The algorithm then iterates over each time step of the planning horizon. Each time step involves first the sampling of a new set of orders $\ExogenousInformation_\CurrentTime$, followed by an enumeration of the feasible order batchings for workers. The state space information and feasible actions are stored as experience, and neural network is utilized to obtain the decision-making model objective coefficients. The task allocation model is then utilized to obtain the optimal solution. Experiences may then be sampled to update the neural network weights, and the network system is updated by implementing the task allocation solution. The algorithm commences by outputting the neural network with the updated weights after simulation.

\section{Experimental Setup} \label{experimentalsetup}
We implement our methods using Python 3.6.13 and run the numerical experiments on Compute Canada Cedar servers~\citep{ComputeCanada}. The linear programming (LP) models are solved using IBM ILOG CPLEX Optimization Studio, version 12.10.0. In what follows, we provide a detailed overview of the experimental settings, which are important for evaluating the effectiveness of our \NeurADP~policy. This includes an in-depth description of the warehouse environment and the data generation, followed by an explanation of the myopic benchmark policies that are used in our comparative analysis.
\subsection{Dataset Description}
\begin{figure}[!ht]
\begin{center}
\includegraphics[width=0.909\textwidth]{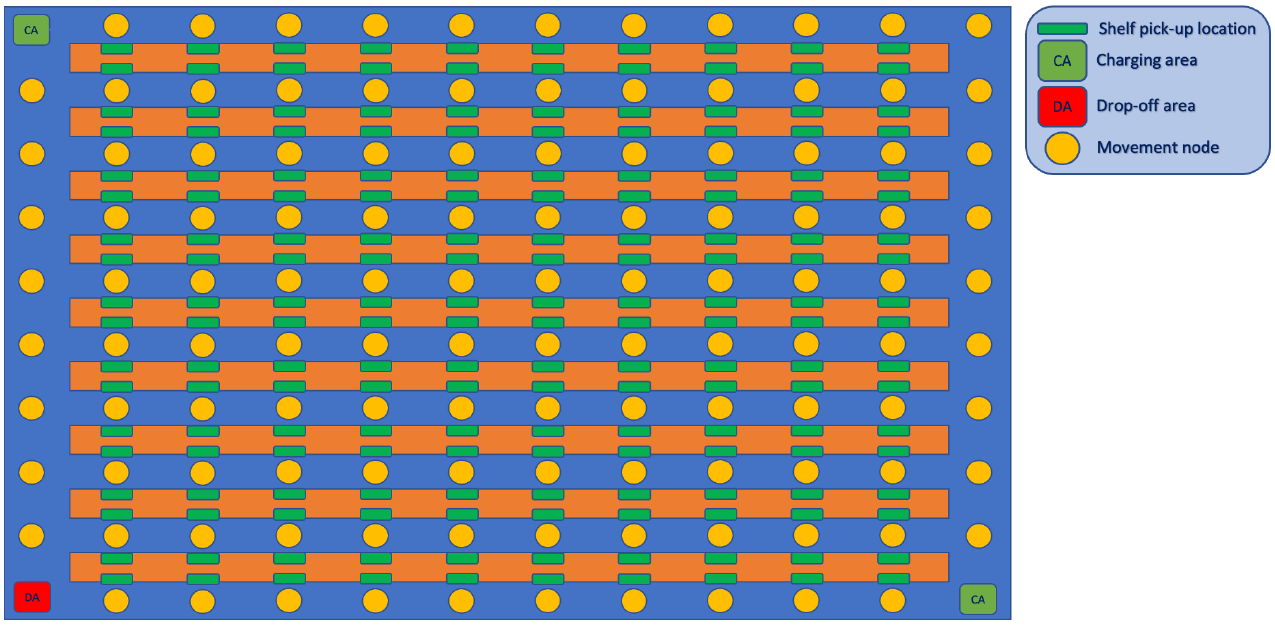}
\end{center}
\caption{A representative warehouse layout.}
\label{fig:Warehouse}
\end{figure}

Our dataset is based on a grid-patterned warehouse which features 9 shelve corridors, each with 20 pick-up locations, amounting to a total of 180 pick-up locations. The warehouse layout includes a single drop-off zone located in the bottom-left corner, along with two charging stations situated in the bottom-right and top-left corners. The warehouse workers and AGVs navigate through aisles spaced between these shelves. In our base-case scenario, the time taken to travel from one movement node to another, including to and from the charging and drop-off points, is uniformly set at 30 seconds. A visual representation of a representative warehouse layout is provided in Figure~\ref{fig:Warehouse}.

\begin{figure}[!ht]
\begin{center}
\includegraphics[width=0.509\textwidth]{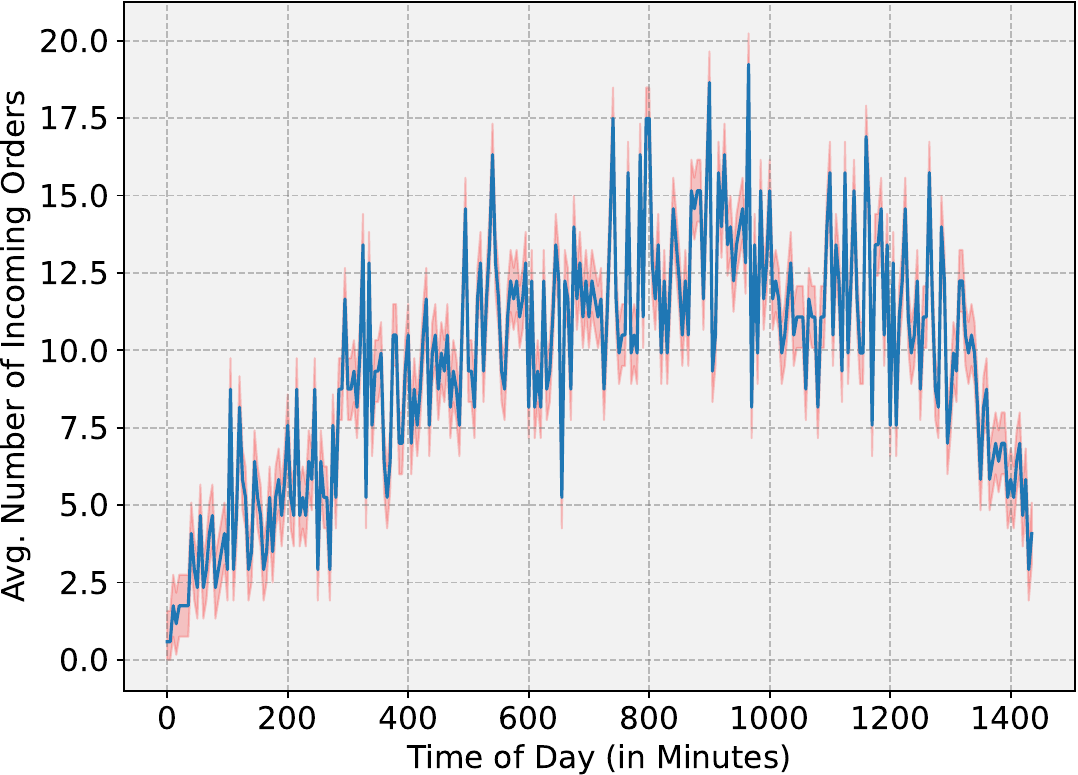}
\end{center}
\caption{Order distribution with 1-order standard deviation band.}
\label{fig:OrderDistribution}
\end{figure}

We use synthetic order arrival data in our numerical experiments.
Specifically, we assume that orders arrive to the system in a stochastic manner, adhering to a left-skewed beta distribution characterized by parameters $\alpha=5$ and $\beta=2$. This distribution, illustrating the average quantity of orders received at each decision epoch, is depicted in Figure~\ref{fig:OrderDistribution}, accompanied by a band representing a standard deviation of one order. Consequently, during each decision epoch, the mean number of arriving orders is utilized as the central value for a normal distribution with a standard deviation of one. The resulting value, rounded to the nearest integer, represents the count of orders arriving at that specific time step for a given simulation iteration. Additionally, the likelihood of an order being requested from any of the 180 pick-up locations at a particular time is determined by sampling from a Poisson distribution with a mean of 1. These probabilities are subsequently normalized to sum to one, transforming them into a proportional distribution. This normalized distribution effectively mirrors the comparative likelihood of order arrivals at each location for every time slot, ensuring a balanced and realistic representation of order frequencies across the network.

\subsection{Benchmark Policies}
Myopic strategies typically prioritize immediate outcomes, neglecting the potential future implications of choices. Such methods are advantageous in complex, rapidly changing environments where developing a comprehensive, optimal strategy is not practical. While myopic policies can be useful for quick decision-making, they might not be ideal for achieving the best long-term results. As a result, they are frequently employed as basic reference points or benchmarks in most cases. Previous studies have predominantly employed simplistic rule-based methodologies for allocating incoming orders among different types of workers. This often involves categorizing orders and then designating these categories to specific worker types \citep{winkelhaus2022hybrid, zhang2023ergonomic, zhang2021evaluation}. 

We consider two distinct sets of myopic policies in our experimental study: \MyopicILP~ and \MyopicHeuristic. Each policy is designed with the primary goal of maximizing the number of orders served in the immediate time step. The \MyopicILP~ policy employs a linear programming approach which adheres to constraints \eqref{m:IPhumanactions}-\eqref{m:variableDefinitions}, ignoring the expected future rewards and only maximizing immediate rewards. Conversely, the \MyopicHeuristic~ approach involves a two-step decision-making process. The first step determines whether orders should be first allocated to human workers or AGVs. The second step involves deciding the optimal time for AGVs to recharge their batteries. Within this framework, \MyopicHeuristic~ policies are further extended: Myopic-HF, which prioritizes assigning orders to humans before AGVs, and Myopic-RF, which does the opposite, giving preference to AGVs for order fulfillment. Additionally, these policies include a battery management component, exemplified by a policy like Myopic-HF-20, which dictates that AGVs not currently assigned to orders should proceed to recharge if their battery levels fall below 20\%. By incorporating both \MyopicILP~and \MyopicHeuristic~approaches in task allocation, we establish a comprehensive baseline in our comparative analysis with our \NeurADP~policy, which allows us to provide an in-depth analysis of the performance of both methods. It also aids in understanding how different strategies for order/task allocation can impact the overall performance of a hybrid workforce in warehouse operations.

We note that the recharging strategies adopted in the two sets of myopic policies are inspired by the charging schemes outlined in existing literature (e.g., see \citep{mchaney1995modelling}). Specifically, the \MyopicILP~policies employ an ``opportunity'' charging approach, allowing AGVs to charge during periods of inactivity. In contrast, the \MyopicHeuristic~policies adopt an ``automatic'' charging strategy, where AGVs are directed to charge only when their battery levels drop below a predefined threshold. Although AGVs in industrial settings typically operate under the latter scheme~\citep{de2020resource}, our exploration of both strategies aims to establish a thorough baseline against which the efficacy of our \NeurADP~policies can be assessed. Hence, this approach allows for a more comprehensive understanding of how different charging strategies impact AGV performance in a warehouse environment.

\section{Results}\label{results}
We evaluate the results from our numerical experiments with respect to five primary inputs: the number of workers, the allowed delay time, the worker capacity, the speed at which each type of worker performs their tasks, and the availability of orders to both human and AGV workers. The number of workers is obtained by accounting for all the workers operating within a 24-hour period, while the delay time represents the maximum duration of time a worker has from an order's entry into the system to their drop-off at the drop-off area. This duration is used to determine the order deadline. 
The worker capacity specifies the maximum number of orders a worker can carry simultaneously, while the speed of workers represents the speed at which they are able to perform their tasks. Finally, the availability of orders specifies the percentage of incoming orders which are able to be handled by both humans and robots. Due to the limitations of AGVs being able to handle items of certain shapes and sizes, we find it important to consider scenarios where only humans are able to handle certain orders, adding complexity to the batching and matching of orders together and to workers. In our \textit{baseline configuration}, we include 5 human workers and 5 AGV workers, set a maximum allowable delay time of 15 minutes, and a maximum worker capacity of 2 orders for both humans and AGVs. Furthermore, we set the travel time of all workers between edges to be 30 seconds and consider all orders to be handleable by both types of workers. Finally, we consider a battery deterioration rate of 0.5\% per minute for AGVs, as well as a charging replenishment rate of 5\% per minute. We begin by examining the baseline configuration to identify the most suitable benchmark policies from the \MyopicILP~and \MyopicHeuristic~policy classes, which are later used in the comparative analysis with the \NeurADP~policy. We then provide the numerical results for alternative configurations.

\subsection{Baseline Configuration}

Our primary benchmark policies exhibit several variations in the matching process between orders and workers. In the \MyopicILP~policy, we utilize an ILP to make decisions on assigning batches of incoming orders with workers, as well as making charging decisions for AGVs. However, in the \MyopicHeuristic~policy, we utilize a heuristic-based mechanism for deciding on which worker to assign orders to first, and then deciding on when to assign unassigned AGVs to recharge their batteries. We consider six variations of the \MyopicHeuristic~policies. Once again, ``HF'' represents the matching of orders with humans first, while ``RF'' represents the matching of orders with AGVs first. Additionally, we consider the battery threshold at which AGVs, which are not currently assigned to orders, are directed to the charging areas to recharge their batteries. We consider the cases where humans are assigned to the orders first, and where robots are assigned to recharge their batteries at thresholds of 20\%, 40\%, and 60\%, as well as the scenarios where robots are first assigned orders, and again where the thresholds are 20\%, 40\%, and 60\%.

Table~\ref{table:base_table} provides the outcomes for the considered policy variants for the baseline configuration, which are obtained as the average statistics when the policies are evaluated on 50 test days. The table includes the ``Orders Seen'' value, which indicates the average number of daily orders seen, as well as the number of orders fulfilled by each policy, denoted as ``Orders Filled'', with a standard deviation provided for each policy. Furthermore, the percentage increase of the \NeurADP~policy compared to the other benchmark policies is included in the right-most column, labeled ``\% Incr. NeurADP''. This metric is calculated by subtracting the average number of orders fulfilled by the benchmark policies from the average number of orders fulfilled by the \NeurADP~policy, then dividing the outcome by the ``Orders Seen'' value, multiplied by 100.

\setlength{\tabcolsep}{1.5pt}
\renewcommand{\arraystretch}{1.15}
\begin{table}[!ht]
\centering
\caption{Performance of different policies for the \textit{baseline configuration} (avg. number of orders fulfilled over 50-day test window is provided as mean$\pm$stdev).}
\label{table:base_table}
\resizebox{0.605\textwidth}{!}{
\begin{tabular}{P{0.165\textwidth}L{0.15\textwidth}L{0.20\textwidth}C{0.23\textwidth}}
\toprule
\textit{\textbf{Policy}} & \textit{\textbf{Orders Seen}} & \textit{\textbf{Orders Filled}} & \textit{\textbf{\% Incr. NeurADP}}\\
\midrule
\textit{\textbf{NeurADP}} & 2,618.25 & 2,051.46 $\pm$ 15.16 & -\\
\midrule
\textit{\textbf{Myopic-ILP}} & - & 1,966.66 $\pm$ 14.82 & +3.24\\
\textit{\textbf{Myopic-HF-20}} & - & 1,946.02 $\pm$ 15.73 & +4.03\\
\textit{\textbf{Myopic-HF-40}} & - & 1,944.92 $\pm$ 16.06 & +4.07\\
\textit{\textbf{Myopic-HF-60}} & - & 1,942.28 $\pm$ 15.16 & +4.17\\
\textit{\textbf{Myopic-RF-20}} & - & 1,918.80 $\pm$ 13.42 & +5.07\\
\textit{\textbf{Myopic-RF-40}} & - & 1,919.08 $\pm$ 13.77 & +5.05\\
\textit{\textbf{Myopic-RF-60}} & - & 1,918.06 $\pm$ 15.72 & +5.09\\
\bottomrule
\end{tabular}
}
\end{table}

We observe that the \NeurADP~policy noticeably outperforms all benchmark policy variants in both the \MyopicILP~and \MyopicHeuristic~scenarios. This superior performance is largely due to its more nuanced approach in effectively pairing order batches with workers and strategically deciding when to assign AGVs for battery recharging, taking into account the downstream effects of its actions. Additionally, the \MyopicILP~policy demonstrates notable superiority over \MyopicHeuristic-based policies. Its use of ILP for simultaneous decision-making in order matching and charging seems more adept for intricate strategies compared to the simplistic assumptions underlying heuristic-based myopic strategies, which do not vary according to scenarios in order matching and charging decisions. Moreover, we find that strategies that give precedence to humans over AGVs for order assignments lead to more favorable results. This effectiveness primarily stems from the reduced frequency of task assignments to AGVs in human-first policies. In scenarios where orders are assigned to the AGVs first, they tend to engage more often in delivery tasks, leading to shorter and less effective charging periods. This is evidenced by the comparison in average AGV battery life: 37.97\% in the Human-First (HF-20) scenario versus 30.83\% in the Robot-First (RF-20) scenario. Moreover, the average count of AGVs charging at any point is notably lower in the HF scenario (1.76) compared to the RF scenario (1.96). This suggests that AGVs in the latter scenario undergo more frequent charging cycles, which are truncated due to commitments to incoming order delivery tasks. Furthermore, we note that among the various benchmark policy variations, \MyopicILP~and Myopic-HF-20 consistently yield the best performance. As a result, we utilize \MyopicILP~and Myopic-HF-20 as our benchmark policies for the remainder of our experiments.

\subsection{Impact of Number of Workers}
We examine the impact of varying the composition of human and AGV workers on the order picking process to better assess the operational effectiveness of our proposed approach. Specifically, we explore three scenarios: all ten workers being humans, a split of five humans and five AGVs, and all ten workers being AGVs. The results, detailed in Table~\ref{table:agents_table}, demonstrate that the \NeurADP~policy consistently outperforms all baseline policies across different worker configurations. 

\setlength{\tabcolsep}{4.5pt}
\renewcommand{\arraystretch}{1.15}
\begin{table}[!ht]
\centering
\caption{Impact of number of workers on order fulfillment (avg. number of orders fulfilled over 50-day test window is provided as mean$\pm$stdev).}
\label{table:agents_table}
\resizebox{0.925\textwidth}{!}{
\begin{tabular}{P{0.135\textwidth}L{0.102\textwidth}L{0.192\textwidth}C{0.18\textwidth}C{0.192\textwidth}L{0.15\textwidth}L{0.162\textwidth}}
\toprule
\textit{\textbf{Number of Workers}} & \textit{\textbf{Orders Seen}} & \textit{\textbf{NeurADP Filled}} & \textit{\textbf{Myopic-ILP Filled}} & \textit{\textbf{Myopic-HF-20 Filled}} & \textit{\textbf{\% Incr. Over Myopic-ILP}} & \textit{\textbf{\% Incr. Over Myopic-HF-20}}\\
\midrule
\textit{\textbf{10 Humans, 0 AGVs}} & 2,618.88 & 2,066.24 $\pm$ 16.38 & 2,019.24 $\pm$ 15.25 & 2,027.26 $\pm$ 14.42 & +1.79 & +1.49\\
\textit{\textbf{5 Humans, 5 AGVs}} & - & 2,051.46 $\pm$ 15.16 & 1,966.66 $\pm$ 14.82 & 1,946.02 $\pm$ 15.73 & +3.24 & +4.03\\
\textit{\textbf{0 Humans, 10 AGVs}} & -
 & 2,007.36 $\pm$ 15.02 & 1,880.76 $\pm$ 16.40 & 1,826.42 $\pm$ 15.44 & +4.83 & +6.91\\
\bottomrule
\end{tabular}
}
\end{table}

\begin{figure}[!ht]
\begin{center}
\includegraphics[width=0.759\textwidth]{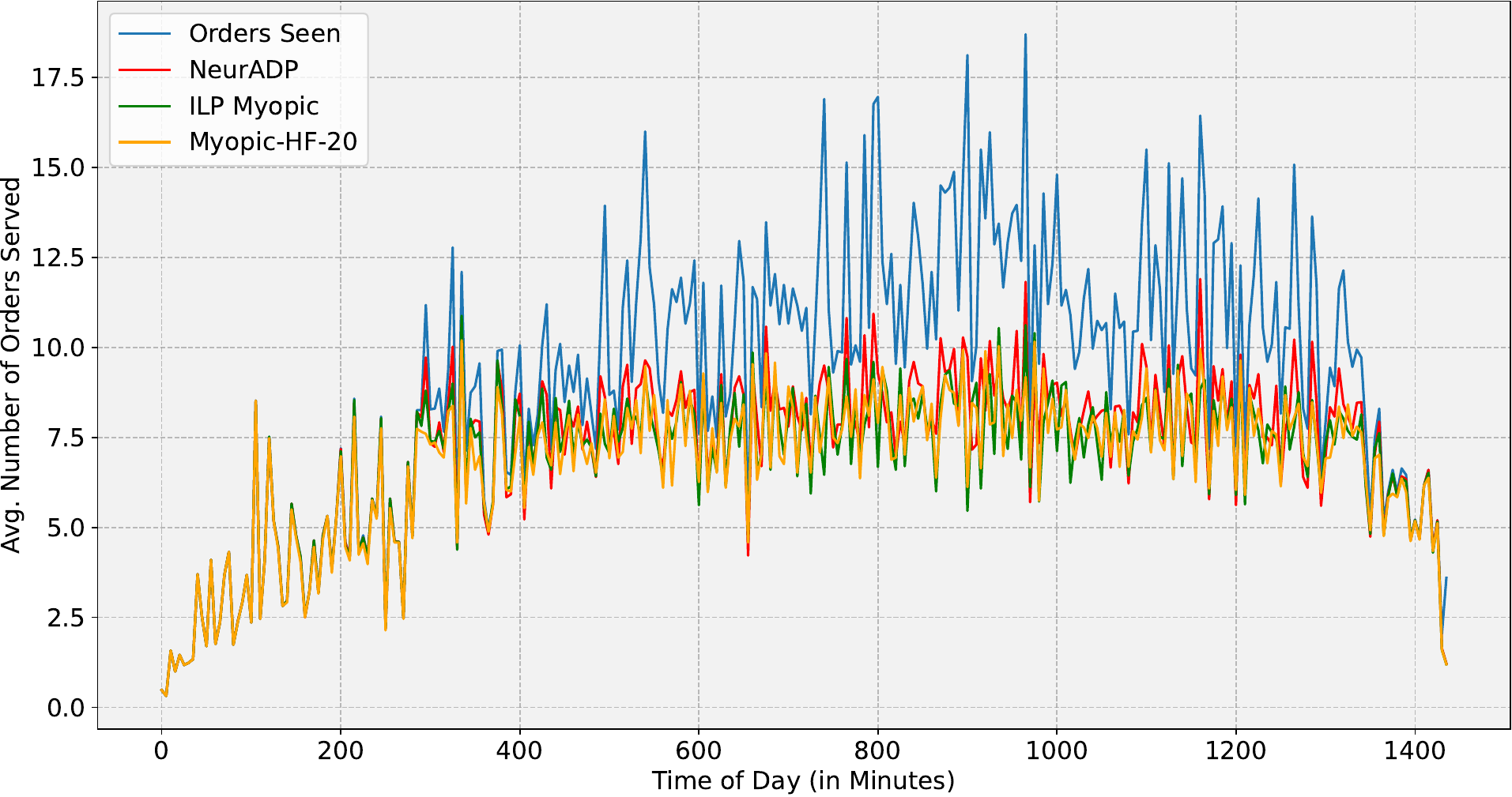}
\end{center}
\caption{Orders fulfilled by each policy in the baseline scenario.}
\label{fig:OrdersServedBase}
\end{figure}

Figure~\ref{fig:OrdersServedBase} presents the number of orders fulfilled in the scenario with an equal split of human and AGV workers. We observe that \NeurADP~maintains superiority over the benchmark policies irrespective of the worker type mix. Notably, there is a trend of decreasing order fulfillment efficiency as the proportion of human workers diminishes. This decline can be attributed to the inherent limitations of AGVs, such as the need for battery charging, which is not a constraint for human workers. Consequently, humans are generally more efficient in handling orders under the given settings. Moreover, there is a noticeable reduction in the performance gap between \NeurADP~and baseline policies as the number of human workers increases. This narrowing of the gap could be linked to the simplification of decision-making processes when AGVs are less involved. In environments dominated by human workers, decision-making becomes less complex, reducing the disparity with simpler myopic policies that may not consider long-term consequences or rely on basic heuristics. In contrast, scenarios including AGVs demand more nuanced decision-making, especially regarding charging strategies, thereby highlighting the advanced decision-making capabilities of \NeurADP~more prominently.

\subsection{Impact of Worker Speed}
We next assess the impact of worker speed on order fulfillment and charging efficiency. Recognizing that AGVs operate at different speeds in various environments, we investigate how these variations in worker efficiency influence the performance of different policies. Specifically, we explore the effects of differing speeds for traversing the warehouse, including scenarios where humans complete an edge in 30 seconds and AGVs in one minute (humans faster than AGVs), both humans and AGVs taking 30 seconds (equal speed), and humans taking one minute while AGVs take 30 seconds (AGVs faster than humans). The results from this experiment are summarized in Table~\ref{table:speed_table}.

\setlength{\tabcolsep}{4.5pt}
\renewcommand{\arraystretch}{1.15}
\begin{table}[!ht]
\centering
\caption{Impact of worker travel time on order fulfillment (avg. number of orders fulfilled over 50-day test window is provided as mean$\pm$stdev).}
\label{table:speed_table}
\resizebox{0.925\textwidth}{!}{
\begin{tabular}{P{0.185\textwidth}L{0.12\textwidth}L{0.182\textwidth}C{0.18\textwidth}C{0.195\textwidth}L{0.15\textwidth}L{0.165\textwidth}}
\toprule
\textit{\textbf{Worker Speed}} & \textit{\textbf{Orders Seen}} & \textit{\textbf{NeurADP Filled}} & \textit{\textbf{Myopic-ILP Filled}} & \textit{\textbf{Myopic-HF-20 Filled}} & \textit{\textbf{\% Incr. Over Myopic-ILP}} & \textit{\textbf{\% Incr. Over Myopic-HF-20}}\\
\midrule
\textit{\textbf{30~Sec.~Human, 1~Min.~AGV}} & 2,618.88 & 1,711.06 $\pm$ 19.02 & 1,651.3 $\pm$ 14.34 & 1,510.86 $\pm$ 16.38 & +2.28 & +7.64\\
\textit{\textbf{30~Sec.~Human, 30~Sec.~AGV}} & - & 2,051.46 $\pm$ 15.16 & 1,966.66 $\pm$ 14.82 & 1,946.02 $\pm$ 15.73 & +3.24 & +4.03\\
\textit{\textbf{1~Min.~Human, 30~Sec.~AGV}} & -
 & 1,617.18 $\pm$ 25.15 & 1,552.20 $\pm$ 32.49 & 1,566.28 $\pm$ 35.11 & +2.48 & +1.94\\
\bottomrule
\end{tabular}
}
\end{table}

\begin{figure}[!ht]
\begin{center}
\includegraphics[width=0.759\textwidth]{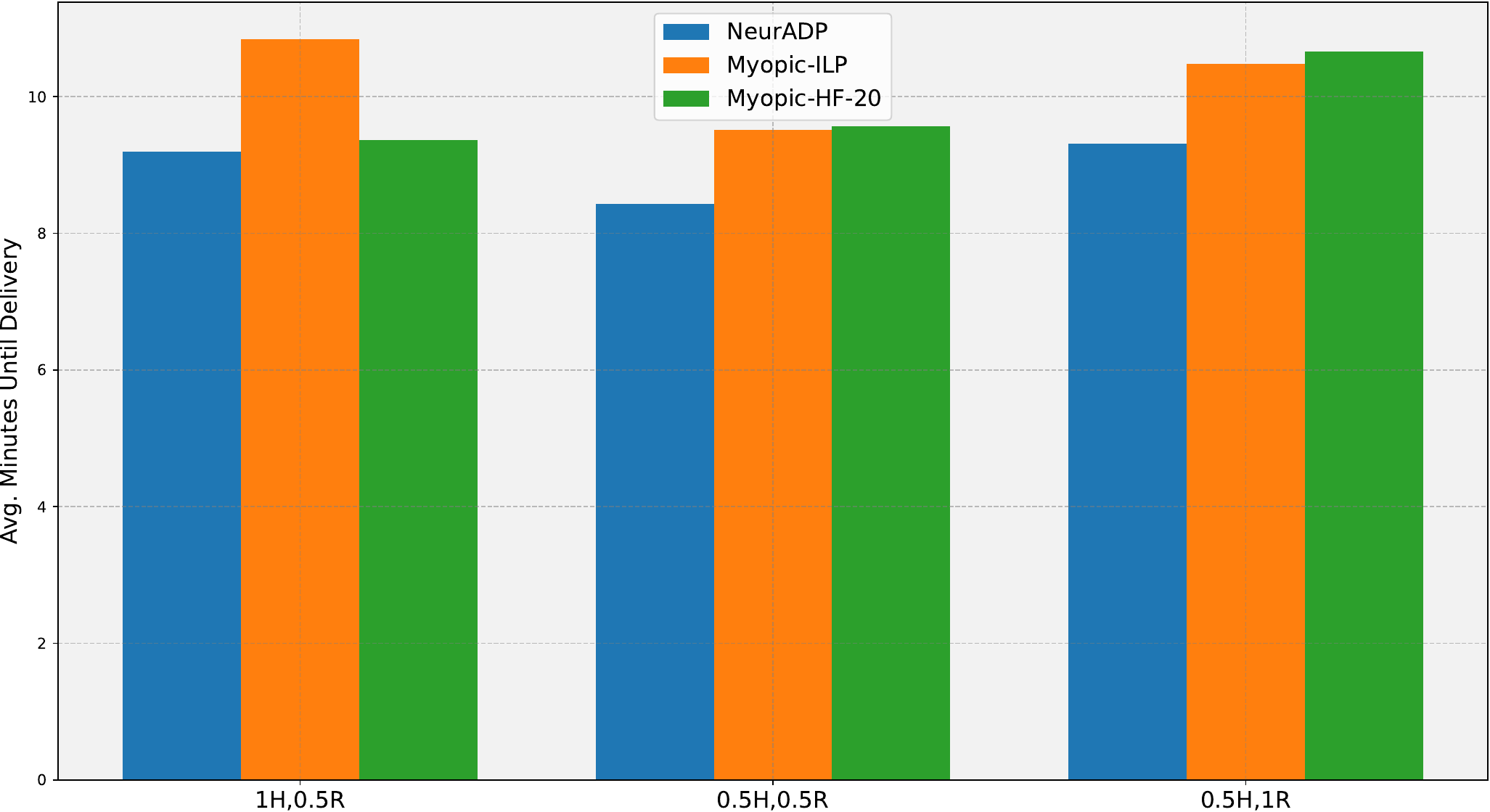}
\end{center}
\caption{Average delivery times of orders for each policy. (``0.5H,1R'' refers 30 second travel time for humans, 1 minute AGVS, ``0.5H,0.5R'' 30 seconds for each, and ``1H,0.5R'' as 1 minute travel time for humans, 30 seconds for AGVs).}
\label{fig:DeliveryTimes}
\end{figure}

Across these scenarios, \NeurADP~consistently outperforms both \MyopicILP~and Myopic-HF-20 policies. This indicates that its neural network-based decision-making adeptly handles the complex dynamics of the warehouse, encompassing variations in worker speed, order batching, and charging schedules. Moreover, \NeurADP~not only serves more orders but also achieves quicker deliveries. Specifically, as depicted in Figure~\ref{fig:DeliveryTimes}, in the three scenarios, \NeurADP~delivers orders in an average of 9.19, 8.43, and 9.31 minutes, compared to \MyopicILP~at 10.84, 9.51, and 10.48 minutes, and Myopic-HF-20 at 9.36, 9.57, and 10.66 minutes, respectively.
We also observe that the policies generally perform better when humans are more efficient than AGVs. This may be attributed to the fact that slower AGVs take longer to reach charging stations, reducing their availability for order assignments. In contrast, slower humans do not encounter this specific issue, leading to a lesser impact on order service rates. Particularly notable is the marked decrease in performance of the Myopic-HF-20 policy when AGVs are slower, possibly due to its threshold-based battery management leading to frequent and inefficient charging trips. Conversely, \NeurADP~seems adept at managing the increased complexity and diversity in worker speeds, thereby making better decisions.

\subsection{Impact of Delay Time}
We also explore how varying the acceptable delay time for order delivery affects the throughput of orders served. For this purpose, we analyze three scenarios: a 10-minute deadline for order drop-off post-entry into the system, an extended deadline of 15 minutes, and a further extension to 20 minutes. We present the results for these experiments in Table~\ref{table:delay_table}. In all these scenarios, \NeurADP~consistently outperforms the benchmark policies regarding allowed delay time. A significant increase in the number of orders served is observed when the delay time is extended from 10 to 15 minutes. This suggests that the initial 10-minute window was overly restrictive, preventing the effective batching and assignment of orders to workers due to time constraints. However, extending the deadline further from 15 to 20 minutes does not yield any notable improvements in order throughput for any of the policies. This plateau in performance improvement may be attributed to other limiting factors on workers, such as capacity constraints, which imply that the additional delay allowance does not translate into the capability to handle more orders or to fulfill orders that were previously unmanageable.

\setlength{\tabcolsep}{4.5pt}
\renewcommand{\arraystretch}{1.15}
\begin{table}[!ht]
\centering
\caption{Impact of delay time on order fulfillment (avg. number of orders fulfilled over 50-day test window is provided as mean$\pm$stdev).}
\label{table:delay_table}
\resizebox{0.925\textwidth}{!}{
\begin{tabular}{P{0.125\textwidth}L{0.12\textwidth}L{0.182\textwidth}C{0.18\textwidth}C{0.195\textwidth}L{0.15\textwidth}L{0.165\textwidth}}
\toprule
\textit{\textbf{Delay Time}} & \textit{\textbf{Orders Seen}} & \textit{\textbf{NeurADP Filled}} & \textit{\textbf{Myopic-ILP Filled}} & \textit{\textbf{Myopic-HF-20 Filled}} & \textit{\textbf{\% Incr. Over Myopic-ILP}} & \textit{\textbf{\% Incr. Over Myopic-HF-20}}\\
\midrule
\textit{\textbf{10 Minutes}} & 2,618.88 & 1,846.62 $\pm$ 21.60 & 1,722.38 $\pm$ 26.23 & 1,695.52 $\pm$ 22.20 & +4.74 & +5.77\\
\textit{\textbf{15 Minutes}} & - & 2,051.46 $\pm$ 15.16 & 1,966.66 $\pm$ 14.82 & 1,946.02 $\pm$ 15.73 & +3.24 & +4.03\\
\textit{\textbf{20 Minutes}} & -
 & 2,053.70 $\pm$ 14.17 & 1,935.16 $\pm$ 16.15 & 1,934.34 $\pm$ 15.58 & +4.53 & +4.56\\
\bottomrule
\end{tabular}
}
\end{table}

\subsection{Impact of Worker Capacity}
The carrying capacity of the workers during pickups can also have a significant impact on the effectiveness of the task allocation policies.
In this experiment, we adjust the capacities of human and AGV workers in the following configurations: AGVs with a capacity of 2, humans with a capacity of 3, both AGVs and humans with a capacity of 2, and finally AGVs with a capacity of 3 and humans with 2. The results from this experiment are summarized in Table~\ref{table:capacity_table}. Initially, we observe that the \NeurADP~policy surpasses benchmark policies across all capacity variations. Further, increasing the carrying capacity of either AGVs or humans boosts the number of orders fulfilled by all policies, indicating that current order servicing is constrained by the carrying capacity of the workers. Moreover, increasing the capacity of human workers appears to be more advantageous for order servicing than doing the same for AGVs. This discrepancy arises because AGVs are limited by the need to periodically recharge their batteries, which can restrict their ability to efficiently utilize additional capacity. Notably, the \NeurADP~policy shows a more pronounced improvement when worker capacities are increased. This improvement can be attributed to the policy's enhanced ability to make more strategic decisions. With fewer constraints on matching order batches with workers, \NeurADP~has greater flexibility to optimize decisions, whereas, in scenarios with stricter capacity limitations, the scope for significant decision-making improvements is more limited.

\setlength{\tabcolsep}{3.0pt}
\renewcommand{\arraystretch}{1.15}
\begin{table}[!ht]
\centering
\caption{Impact of worker capacity on order fulfillment (avg. number of orders fulfilled over 50-day test window is provided as mean$\pm$stdev).}
\label{table:capacity_table}
\resizebox{0.925\textwidth}{!}{
\begin{tabular}{P{0.115\textwidth}L{0.12\textwidth}L{0.182\textwidth}C{0.18\textwidth}C{0.195\textwidth}L{0.15\textwidth}L{0.165\textwidth}}
\toprule
\textit{\textbf{Worker Capacity}} & \textit{\textbf{Orders Seen}} & \textit{\textbf{NeurADP Filled}} & \textit{\textbf{Myopic-ILP Filled}} & \textit{\textbf{Myopic-HF-20 Filled}} & \textit{\textbf{\% Incr. Over Myopic-ILP}} & \textit{\textbf{\% Incr. Over Myopic-HF-20}}\\
\midrule
\textit{\textbf{AGV:~2, Human:~3}} & 2,618.88 & 2,284.18 $\pm$ 14.10 & 2,140.88 $\pm$ 14.58 & 2,125.1 $\pm$ 15.13 & +5.47 & +6.07\\
\textit{\textbf{AGV:~2, Human:~2}} & - & 2,051.46 $\pm$ 15.16 & 1,966.66 $\pm$ 14.82 & 1,946.02 $\pm$ 15.73 & +3.24 & +4.03\\
\textit{\textbf{AGV:~3, Human:~2}} & -
 & 2,262.90 $\pm$ 12.90 & 2,115.02 $\pm$ 15.74 & 2,121.9 $\pm$ 13.71 & +5.65 & +5.38\\
\bottomrule
\end{tabular}
}
\end{table}

\subsection{Impact of Order Availability}
Lastly, we examine how the availability of orders for different worker types affects their ability to fulfill orders. Considering the limited capacity of certain AGVs to handle specific items or orders, we evaluate scenarios where AGVs can only manage a certain percentage of incoming orders. Specifically, we assess cases where 0\%, 20\%, and 40\% of orders can exclusively be handled by humans, with results detailed in Table~\ref{table:availability_table}. This limitation challenges AGVs in processing order batches, as they can only be assigned orders within their handling capacity. Consequently, we observe a decrease in the number of orders served across all policies as the proportion of human-exclusive orders increases. However, the \NeurADP~policy, thanks to its more complex decision-making capabilities, consistently outperforms the baseline policies, especially in optimally utilizing AGVs under these constraints. This is particularly evident in the 40\% scenario, where the humans operating under the \NeurADP~policy serve a comparable number of orders to those in benchmark policies, but the AGVs in the \NeurADP~system manage a significantly higher number of orders. This leads to an overall better performance relative to the benchmark policies. For instance, in the \NeurADP~policy under the 40\% scenario, AGVs serve an average of 866.98 orders, markedly more than the 798.3 and 624.94 orders served under the \MyopicILP~and Myopic-HF-20 policies, respectively. This demonstrates \NeurADP's proficiency in leveraging AGVs effectively, even when faced with order availability restrictions.

\setlength{\tabcolsep}{4.5pt}
\renewcommand{\arraystretch}{1.15}
\begin{table}[!ht]
\centering
\caption{Impact of order availability on order fulfillment (avg. number of orders fulfilled over 50-day test window is provided as mean$\pm$stdev).}
\label{table:availability_table}
\resizebox{0.925\textwidth}{!}{
\begin{tabular}{P{0.145\textwidth}L{0.12\textwidth}L{0.182\textwidth}C{0.18\textwidth}C{0.195\textwidth}L{0.15\textwidth}L{0.165\textwidth}}
\toprule
\textit{\textbf{Order Availability}} & \textit{\textbf{Orders Seen}} & \textit{\textbf{NeurADP Filled}} & \textit{\textbf{Myopic-ILP Filled}} & \textit{\textbf{Myopic-HF-20 Filled}} & \textit{\textbf{\% Incr. Over Myopic-ILP}} & \textit{\textbf{\% Incr. Over Myopic-HF-20}}\\
\midrule
\textit{\textbf{0\%}} & 2,618.88 & 2,051.46 $\pm$ 15.16 & 1,966.66 $\pm$ 14.82 & 1,946.02 $\pm$ 15.73 & +3.24 & +4.03\\
\textit{\textbf{20\%}} & - & 2,038.86 $\pm$ 15.39 & 1,960.06 $\pm$ 15.84 & 1,899.70 $\pm$ 14.85 & +3.01 & +5.31\\
\textit{\textbf{40\%}} & -
 & 2,019.36 $\pm$ 14.05 & 1,899.98 $\pm$ 17.41 & 1,872.00 $\pm$ 14.34 & +4.56 & +5.63\\
\bottomrule
\end{tabular}
}
\end{table}

\section{Conclusion}\label{conclusion}
In this study, we investigate the AGV integration within picker-to-parts warehouse systems, and we model a dynamic warehouse environment where human workers and AGVs work in tandem. For the considered problem setting, we develop a NeurADP-based solution approach, which enables non-myopic decision-making in order allocation and battery management for AGVs. The detailed numerical study underscores the NeurADP approach's efficacy, yielding a significant improvement over traditional myopic and heuristic-based methods. Specifically, the NeurADP policies exhibit superior performance in fulfilling a higher number of orders and enhancing the efficiency of executing order picking tasks, which is particularly evident under varying problem parameters such as worker speed, delay time, and order availability. Furthermore, our analysis provides valuable managerial insights into operating a hybrid warehouse and the importance of intelligent decision-making in a mixed workforce environment.

There exist several study limitations that can be remedied in future research. First, our experiments are conducted with a synthetic dataset that is generated based on real-life warehouse configurations and order arrival patterns. However, extending our numerical study to alternative warehouse configurations and order arrival patterns can help further validate our proposed methods. Moreover, the complexity of human and AGV interaction is modeled under simplifying assumptions that may not capture the full range of behaviors in a real-world setting.
Exploring the psychological and behavioral aspects of human-AGV interaction may provide deeper insights into optimizing the collaborative workspace for both efficiency and worker satisfaction. In terms of modeling and methodological extensions, the learning capacity of the NeurADP can be potentially enhanced for the AGV task allocation problem by considering more complex NN designs and other algorithmic enhancements such as extending individual worker-based value function approximation to the worker sets. Another valuable future research direction in this regard is to extend our proposed approach to solve the integrated AGV management problem that also involves localization, motion planning and path planning steps. Ultimately, as the landscape of warehouse automation continues to evolve, the pursuit of innovative solutions such as NeurADP will remain pivotal in driving the industry forward.

\section*{Disclosure statement}
No potential conflict of interest was reported by the authors.



\singlespacing
\bibliographystyle{elsarticle-harv} 
\bibliography{main_elsevier}

\end{document}